\documentclass[11pt]{article}

\usepackage[a4paper, margin=1in]{geometry}
\usepackage{amsfonts,amsmath,amssymb,amsthm,graphicx,fixmath,latexsym,color}
\usepackage{xcolor}
\usepackage{hyperref,epsfig}
\usepackage{algorithmic}
\usepackage{algorithm}
\usepackage{xspace}

\usepackage{graphicx}
\usepackage{caption,subcaption}

\providecommand{\remove}[1]{}

\newtheorem{theorem}{Theorem}[section]
\newtheorem{lemma}[theorem]{Lemma}

\newtheorem{conjecture}[theorem]{Conjecture}

\newtheorem*{theorem*}{Theorem}
\newtheorem*{lemma*}{Lemma}
\newtheorem*{proposition*}{Proposition}
\newtheorem{observation}[theorem]{Observation}
\newtheorem*{example*}{Example}

\newcommand{\bbox}{\vrule height7pt width4pt depth1pt}
\newcommand{\dsum}{\dot{+}}

\begin{document}
\title{On sets of points in general position that lie on a cubic curve
in the plane and determine lines that can be pierced by few points}

\author{Mehdi Makhul\thanks{Johann Radon Institute for Computational and Applied Mathematics, Austrian Academy of Sciences, Altenberger Str. 69, 4040 Linz, Austria. \texttt{mehdi.makhul@oeaw.ac.at}. Research partially supported by Grant Artin Approximation, Arc-R{\"a}ume, Aufl{\H o}sung von Singularit{\"a}ten FWF P-31336, and the Austrian Science Fund FWF Project P 30405.}
\and
Rom Pinchasi\thanks{Mathematics Department, Technion -- Israel Institute of Technology, Haifa 32000, Israel. \texttt{room@tx.technion.ac.il}. Research partially supported by Grant 409/16 from the Israel Science Foundation. The second author acknowledge the financial support from the Ministry of Educational and Science of the Russian Federation in the framework of MegaGrant no 075-15-2019-1926.}
}

\date{}
\maketitle

\begin{abstract}
Let $P$ be a set of $n$ points in general position in the plane.
Let $R$ be a set of points disjoint from $P$ such that
for every $x,y \in P$ the line through $x$ and $y$ contains a point in $R$.
We show that if $|R| < \frac{3}{2}n$ and $P \cup R$ is contained in a cubic
curve $c$ in the plane, then $P$ has a special property with respect to 
the natural group action on $c$. That is,
$P$ is contained in a coset of a subgroup $H$ of $c$ of cardinality at most
$|R|$. 

We use the same approach to show a similar result in the case where each of
$B$ and $G$ is a set of $n$ points in general position in the plane and 
every line through a point in $B$ and a point in $G$ passes through a point
in $R$. This provides a partial answer to a problem of Karasev.

The bound $|R| < \frac{3}{2}n$ is best possible at least for part of our 
results. Our extremal constructions provide
a counterexample to an old conjecture attributed to Jamison about 
point sets that determine few directions. Jamison conjectured that if $P$ is a set of $n$ points in general position in the plane that determines at most 
$2n-c$ distinct directions, then $P$ is contained in an affine image of
the set of vertices of a regular $m$-gon. 
This conjecture of Jamison is strongly
related to our results in the case the cubic curve $c$ is reducible and 
our results can be used to prove Jamison's conjecture at least when
$m-n$ is in the order of magnitude of $O(\sqrt{n})$.    
\end{abstract}


\section{Introduction}
\label{sec:introduction}

\begin{figure}[ht]
   \centering
    \includegraphics[width=7cm]{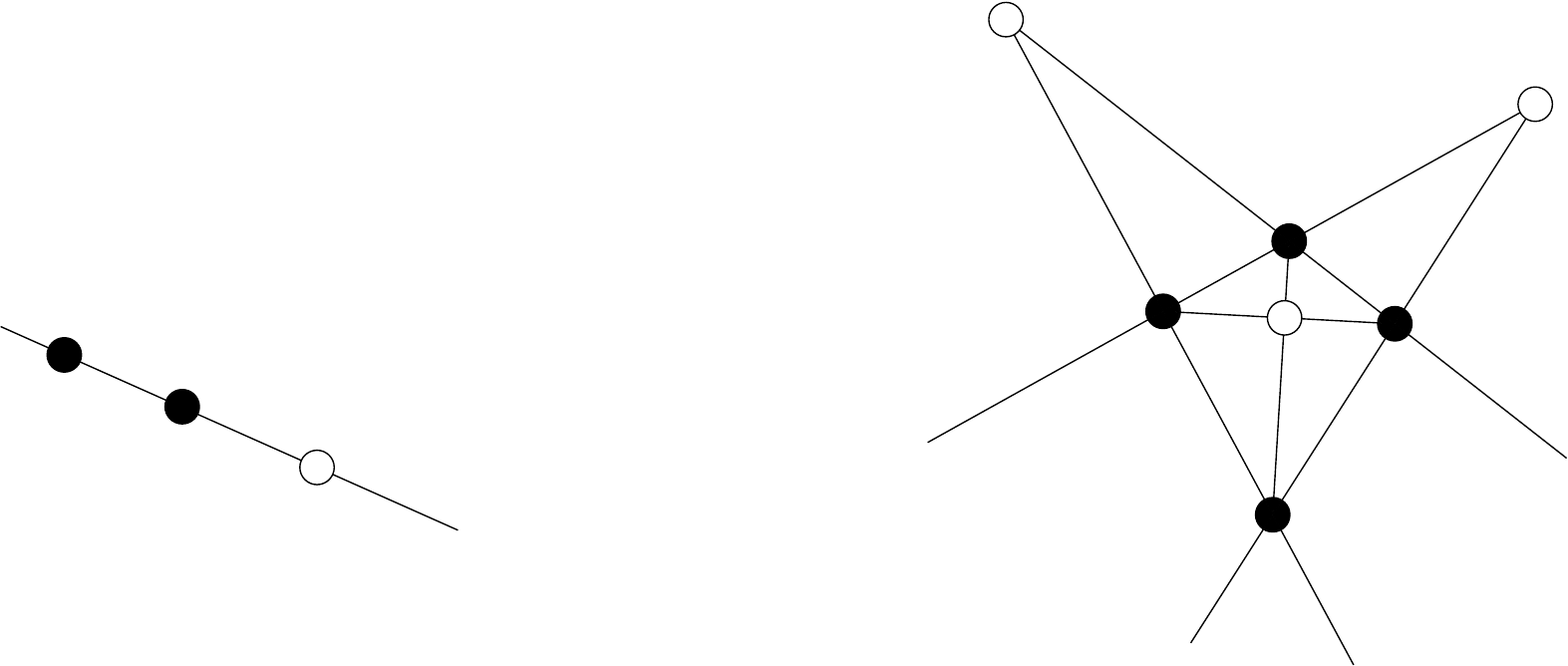}
        \caption{Constructions with $|R|=n-1$ for $n=2,4$.
The points in $P$ are colored black while the points in $R$ are colored white.}
        \label{figure:fig1}
\end{figure}

In \cite{EP78} Erd\H os and Purdy considered the following problem (Problem 3
in \cite{EP78}).
Let $P$ be a set of $n$ points in the plane that lie  
in \emph{general position} in the sense that no three points of $P$
are collinear. Let $R$ be another set of points disjoint from $P$ such that
every line through two points of $P$ contains a point in $R$. Give a lower
bound on $|R|$ in terms of $n$.

If $n$ is odd it is easy to prove the tight bound $|R| \geq n$.
This is because every point in $R$ may be incident to at most
$\frac{n-1}{2}$ of the ${n \choose 2}$ lines determined by $P$.
To observe that this bound is tight let $P$ be the
set of vertices
of a regular $n$-gon and let $R$ be the set of $n$ points on the line
at infinity that correspond to the directions of the edges (and diagonals)
of $P$. This construction is valid also when $n$ is even.

If $n$ is even, a trivial counting argument shows that $|R|$ must be at least
$n-1$. This is because every point in $R$ may be incident to at most
$n/2$ lines determined by $P$. This trivial lower bound for $|R|$ is in fact
sharp in the cases $n=2$ and $n=4$, as can be seen in Figure \ref{figure:fig1}.
Is the bound $|R| \geq n-1$ sharp also for larger values of $n$?

The following theorem proves a conjecture attributed to Erd\H os and Purdy \cite{EP78}. We note that the problem is expilicly stated in \cite{EP78} but the
conjectured lower bound is mistakenly constant times $n$ rather than just $n$.

\begin{theorem}[\cite{ABKPR08,M18,PP20}]\label{conjecture:EP}
Let $P$ be a set of $n$ points in general position in the plane, 
where $n>4$ is even. Assume $R$ is another set of points disjoint from
$P$ such that every line through two points of $P$ contains a point from $R$.
Then $|R| \geq n$.
\end{theorem}

Theorem \ref{conjecture:EP} was first proved in
\cite{ABKPR08} (see Theorem 8 there), as a special case of the solution of the
Magic Configurations conjecture of Murty \cite{Murty71}. The proof
in \cite{ABKPR08} contains a topological argument based on Euler's formula
for planar maps and the discharging method.
An elementary (and long) proof of Theorem
\ref{conjecture:EP} was given by Mili\'{c}evi\'{c} in \cite{M18}.
Probably the ``book proof'' of the Theorem \ref{conjecture:EP} can be found in
\cite{PP20}.

Theorem \ref{conjecture:EP} was proved also over $\mathbb{F}_{p}$
by Blokhuis, Marino, and Mazzocca \cite{BMM14}.

As we have seen, there are constructions of sets $P$ of $n$ points in general
position and sets $R$ of $n$ points not in $P$, such that every line
determined by $P$ passes through a point in $R$. One major question that
arises here is to characterize those sets $P$ in general position for which
there exists a set $R$ with $|R|=|P|$ such that every line that is determined
by $P$ passes through a point in $R$.

\medskip

The following conjecture of Mili\'{c}evi\'{c} in \cite{M18}, 
came up in connection to the
above mentioned Theorem \ref{conjecture:EP}.

\begin{conjecture}\label{conjecture:main}
Let $P$ be a set of $n$ points in general position and let $R$
be a set of $n$ points disjoint from $P$. If every line determined by
$P$ passes through a point in $R$, then $P \cup R$ is contained in
a cubic curve.
\end{conjecture}

A special case of Conjecture \ref{conjecture:main} is proved in \cite{KP20},
where Conjecture \ref{conjecture:main} is proved 
under additional assumption:

\begin{theorem}\label{theorem:KP}
Suppose $P$ is a set of $n$ points in general position in the plane and
$R$ is another set of $n$ points, disjoint from $P$.
If for every $x,y \in P$ there is a point $r \in R$ on the line through $x$
and $y$ and outside the interval determined by $x$ and $y$, then
$P \cup R$ is contained in a cubic curve.
\end{theorem}

We recall that given an irreducible 
cubic curve $c$ in the plane, there is a natural
abelian group structure on $c$ (we refer the reader to \cite{B06} and the references therein). 
In this group structure the sum of three collinear
points on $c$ is equal to~$0$. Moreover, if there is a line $\ell$ that
crosses $c$ at a point $A$ and is tangent to $c$ at a point $B$, then 
$A+B+B=0$.   

In this paper we show that if indeed $P \cup R$ in Conjecture 
\ref{conjecture:main} is contained in a cubic curve $c$ in the plane and
if $|P|$ is not too small, then 
there is a subgroup $H$ of $c$ such that both $P$ and $R$ are cosets 
of $H$. In fact, we can extend it as follows:



\begin{theorem}\label{theorem:main}
Let $P$ be a set of $n>6$ points in general position in the plane.
Let $R$ be another set of less than $\frac{3}{2}n$ points, 
disjoint from $P$ such that 
any line through two points of $P$ passes through a point in $R$.
Assume that $P \cup R$ is contained in an irreducible 
cubic curve $c$ in the plane.
Then there is a subgroup $H$ of $c$ of size at most $|R|$ 
such that $P$ is contained in
a coset of $H$. If $|R|=n$, then both $P$ and $R$ are equal to cosets
of $H$. 
\end{theorem}

It is not hard to consider the situation in 
Theorem \ref{theorem:main} 
also in the case where the cubic curve $c$ is 
reducible. In such a case $c$ is either a union of three lines, or a union
of a quadric and a line. The former case is literally impossible as we assume
$P$ is in general position and large enough. The following easy theorem settles
the case where $c$ is a union of a quadric and a line:

\begin{theorem}\label{theorem:main_reducible}
Let $P$ be a set of $n$ points in general position in the plane
and assume that $n$ is large enough ($n>100$ will work here).
Let $R$ be another set disjoint from $P$ such that 
any line through two points of $P$ passes through a point in~$R$.
Assume that $P \cup R$ is contained in a reducible 
cubic curve $c$ that is a union of a quadric~$Q$ and a line $\ell$.
If $n \leq |R| < \frac{3}{2}n$ and $\ell$ is the line at infinity, 
then $Q$ must be an ellipse, 
$P \subset Q$, and up to an affine transformation~$P$ is a subset
of the set of vertices of a regular $m$-gon for some $m \leq |R|$.
The bound $\frac{3}{2}n$ in the statement of the theorem is the best possible.
\end{theorem}

We remark that for values of $n$ smaller than or equal to $6$ one can indeed 
find some sporadic examples that are contained in a union of three lines. 
For instance, the examples in Figure \ref{figure:fig1} are of 
sets of points contained in a union of at most three lines and satisfy in some 
respect the conditions in Theorem \ref{theorem:main}. 

Another interesting remark is that the bound $|R| < \frac{3}{2}n$ in Theorem 
\ref{theorem:main_reducible} cannot be improved, not even by one unit.
We now present a simple construction showing  this. Before presenting
the construction, the following two very easy observations will be useful.

\begin{observation}\label{observation:regular}
Let $S$ be the set of vertices of a regular $m$-gon.
Then $S$ determines lines that appear in precisely $m$ distinct directions.
\end{observation}

\begin{figure}[ht]
   \centering
    \includegraphics[width=7cm]{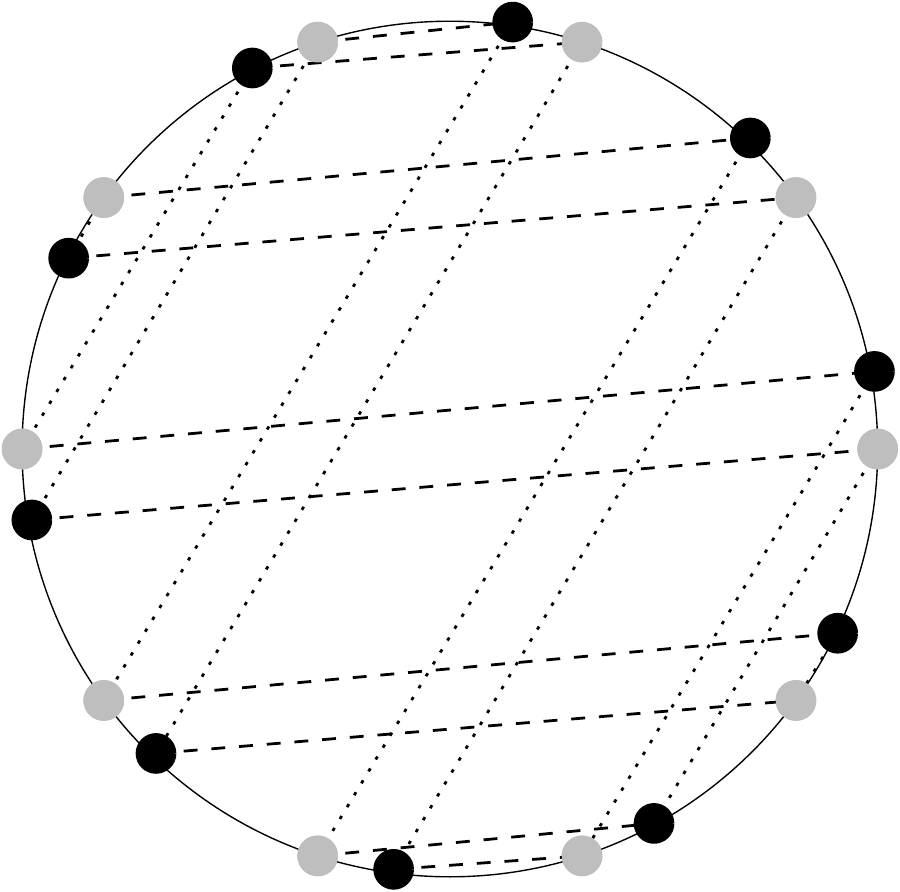}
        \caption{Illustration of Observation \ref{observation:directions} 
for the case $m=10$. The points in $P_{1}$ are colored black while the points in $P_{2}$ are colored gray. The diagonals connecting points of $P_{1}$ and points of $P_{2}$ can be partitioned into $m$ families of
pairwise parallel diagonals.}
        \label{figure:observation}
\end{figure}

\begin{observation}\label{observation:directions}
Let $Q$ be a circle and let $P_{1}$ be the set of vertices of a regular 
$m$-gon inscribed in $Q$. Let $P_{2}$ the set of vertices of another 
regular $m$-gon inscribed in $Q$ such that $P_{1}$ is disjoint from $P_{2}$. 
Then the lines connecting a point from
$P_{1}$ to a point from $P_{2}$ may appear in one of at most $m$ distinct 
directions (see Figure \ref{figure:observation}). 
\end{observation}

We leave it to the reader to verify the validity of 
Observation \ref{observation:regular} and 
Observation \ref{observation:directions}.

In order to see that the bound $|R| < \frac{3}{2}n$ in Theorem 
\ref{theorem:main_reducible} is best possible consider a circle $Q$ and let 
$P_{1}$ be the set of vertices of 
a regular $\frac{n}{2}$-gon inscribed in $Q$ (we assume $n$ is even).
Let $P_{2}$ be a generic rotation (for example, by an angle that is an irrational multiple of $\pi$) 
of $P_{1}$ about the center of $Q$. 
Let $P$ be the union $P_{1} \cup P_{2}$. Therefore, $P$ is a set of $n$
points contained in $Q$. 
By Observation \ref{observation:regular}, the lines passing through two points 
of $P_{1}$ have precisely $\frac{n}{2}$ distinct directions, and
the same is true for the lines passing through two points of $P_{2}$.

From 
Observation \ref{observation:directions} we know that the lines passing through a point of $P_{1}$ and a point of $P_{2}$
have precisely $\frac{n}{2}$ distinct directions. 
We conclude that the lines passing through two points of $P$
have at most and in fact precisely $\frac{3}{2}n$ distinct directions. 
Define $R$ to be the set of $\frac{3}{2}n$ points on the line at 
infinity that correspond to the $\frac{3}{2}n$ directions of lines determined
by the points in $P$. Then every line through two points of $P$
passes through a point of $R$. Moreover, $P$ is not a subset of the 
set of vertices of any regular polygon. 

It is extremely interesting to remark about the relation between Theorem
\ref{theorem:main} and Theorem \ref{theorem:main_reducible} and an old
conjecture attributed to Jamison in \cite{J86}. By a well known theorem 
of Ungar (\cite{Un82}), proving a conjecture of Scott (\cite{Sc70}), every 
set of $n$ points in the plane that is not collinear must determine lines
appearing in at least $2\lfloor \frac{n}{2} \rfloor$ distinct directions.
In \cite{J86} Jamison studied the case where the set $P$ 
of $n$ points is in general
position in the sense that no three of its points are collinear.
In this case it is very easy to prove that $P$ determines at least
$n$ distinct directions. Jamison showed that in the extremal case
if $P$ is a set of $n$ points in general position that determines lines in 
precisely $n$ distinct directions, then $P$ is, up to an affine transformation,
the set of vertices of a regular $n$-gon.
Jamison mentioned in \cite{J86} that is it believed that if a set $P$ of
$n$ points in general position determines $m \leq 2n-c$ distinct directions
(for some large enough absolute constant $c$), then up to an affine
transformation $P$ is contained in the set of vertices of a regular
$m$-gon. This conjecture mentioned in \cite{J86}
was recently addressed in \cite{P18}, where the conjecture
is proved for the case $m=n+1$.

Having shown that the bound $|R| < \frac{3}{2}n$ in 
Theorem \ref{theorem:main_reducible} is best possible we constructed
a set $P$ of $n$ 
points (where we assumed $n$ is even) that determines precisely $\frac{3}{2}n$ 
distinct directions. This set $P$ that we constructed was contained in a circle
but was not a subset of any set of vertices of a regular $m$-gon. This
shows that Jamison's conjecture is false for $m \geq \frac{3}{2}n$ and
one may hope to prove it only for smaller values of $m$.

One possible way to approach Jamison's conjecture is to show that if $P$
is a set of $n$ points that determines less than $\frac{3}{2}n$ distinct 
directions, then $P$ is contained in a conic. 
Then we can take $R$ to be the set of points on the line at infinity that
correspond to the directions determined by $P$ and 
Theorem \ref{theorem:main_reducible} will automatically imply that 
$P$ must be contained (up to an affine transformation) in the set of vertices
of a regular $m$-gon for some $m < \frac{3}{2}n$. 

Such an approach is carried out in an ongoing work 
by the second author and Alexandr Polyanskii 
in the case where $P$ determines at most $n + O(\sqrt{n})$
distinct directions. One can show that such a set $P$ must be 
contained in a conic. Consequently, together with Theorem \ref{theorem:main_reducible}, this extends significantly the result of C. Pilatte in \cite{P18} 
and in fact proves Jamison's conjecture in these cases. Very recently,
independently, Pillate observed a similar improvement in the same spirit
in \cite{P18b}.

\bigskip
Next we turn to present the second part of our paper that is of similar nature.
A bipartite version of the situation that appears in 
Conjecture \ref{conjecture:main} was raised by Karasev:
Assume $P$ is a set of $2n$ points that is the union of a set $B$ of
$n$ blue and a set $G$ of $n$ green points. The set $R$ is a set of $n$ red
points and we assume that the sets $B, G$, and $R$ are pairwise disjoint.
We also assume that the set
$P=B \cup G$ is in general position in the sense that no
three of its points are collinear. Assume that every line through a point in $B$
and a point in $G$ contains also a point from $R$. The problem of characterizing
the set $B \cup G \cup R$ was raised by Karasev (in a dual version) 
in \cite{K18} (see Problem 6.2 there).

In this context the following 
analogous conjecture to Conjecture \ref{conjecture:main} 
appears in \cite{KP20}:

\begin{conjecture}[\cite{KP20}]\label{conjecture:Karasev}
Let $B$, $G$, and $R$ be three sets of points in the plane, each of which
is in general position and has size $n$. Assume that every line passing through
two points from two different sets passes also through a point from the
third set. Then $B \cup G \cup R$ lies on a cubic curve.
\end{conjecture}

Conjecture \ref{conjecture:Karasev} is proved in \cite{KP20} under a similar
additional assumption as in Theorem \ref{theorem:KP}.

\begin{theorem}\cite{KP20}\label{theorem:bipartite}
Let $B, G$, and $R$ be three pairwise disjoint sets of points in the
plane. Assume that $B \cup G$ is in general position and
$|B|=|G|=|R|=n$. If every line through a point $b \in B$ and a point
$g \in G$ contains a point $r \in R$ that does not lie between $b$ and $g$,
then $B \cup G \cup R$ lies on some cubic curve.
\end{theorem}

Similar to Theorem \ref{theorem:main}, we can prove something about the
algebraic structure of $B \cup G \cup R$ in Conjecture \ref{conjecture:Karasev}
and in particular in Theorem \ref{theorem:bipartite}.


\begin{theorem}\label{theorem:mainb}
Let $B, G$, and $R$ be three pairwise disjoint sets of points in the
plane such that~$B \cup G$ is in general position and 
every line through a point in $B$ and a point in $G$ passes
through a point in $R$. 
Assume $|B|=|G|=n$, $|R| < \frac{3}{2}n$ and $B \cup G \cup R$ is 
contained in an 
irreducible cubic curve $c$,
then the sets $B, G$ are contained in cosets of the same subgroup
$H$ of~$c$ of cardinality at most $|R|$. 
\end{theorem}

We can consider Theorem \ref{theorem:mainb} also in the case $c$ is a reducible
cubic curve. Here too if $c$ is a union of three lines and $|B|$ and $|G|$ are 
greater than $6$ 
we get a contradiction to the assumption that both $B$ and $G$ are in 
general position.

\begin{figure}[ht]
   \centering
    \includegraphics[width=7cm]{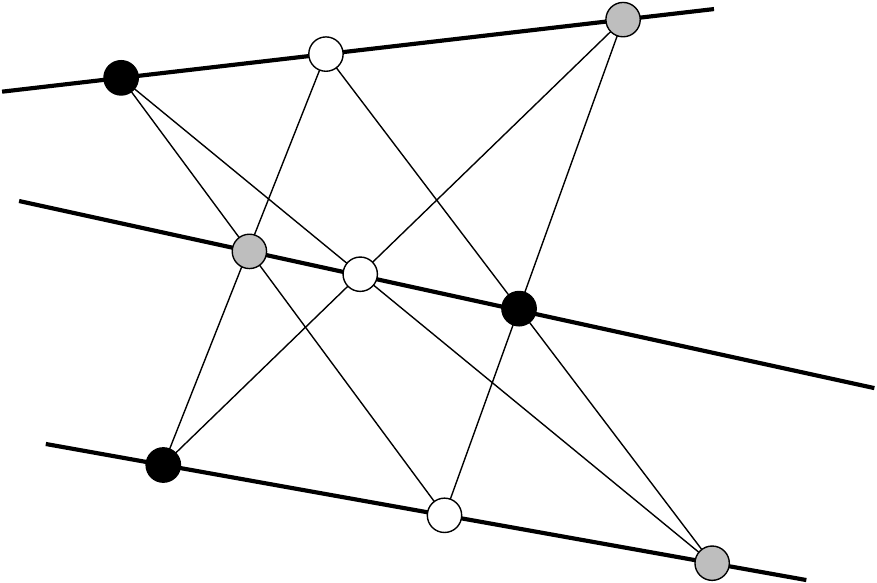}
        \caption{Pappus' Theorem gives rise to a small example for
Theorem \ref{theorem:mainb} that is contained in three lines.
The points of $B$ are colored black, the points of $G$ are colored gray,
and the points of $R$ are colored white.}
        \label{figure:fig2}
\end{figure}

It is interesting to remark that when $|B|=|G| \leq 6$ we may get 
examples satisfying 
the conditions in Theorem \ref{theorem:mainb}. One easy example is
the case $n=1$ of three collinear points 
(one of $B$, one of $G$, and one of $R$). A more interesting example is
Pappus' Theorem that is illustrated in Figure \ref{figure:fig2}. 
In this example
each of the sets $G, B$, and $R$ consists of $3$ points and they satisfy
the conditions in Theorem \ref{theorem:mainb}. The union $B \cup G \cup R$
is contained in a union of three lines.
 
In the case $n>6$, $B \cup G$ cannot be contained in a union of 
three lines and the only possibility for a reducible cubic curve
containing $B \cup G$ is where $c$ is the 
union of a quadric and a line. This case is studied in the following theorem.

\begin{theorem}\label{theorem:mainb_reducible}
Let $B, G$, and $R$ be three pairwise disjoint sets of points in the
plane such that~$B \cup G$ is in general position and 
every line through a point in $B$ and a point in $G$ passes
through a point in $R$. 
Assume $|B|=|G|=n>6$, $|R| < \frac{3}{2}n$, and 
$B \cup G \cup R$ is contained in a reducible 
cubic curve $c$ that is a union of a quadric $Q$ and a line $\ell$,  
then the sets $B, G \subset Q$, and $Q$ must be an ellipse. Moreover,
if $\ell$ is the line at infinity, then 
up to an affine transformation that takes $Q$ to a circle, 
each of $B$ and $G$ is a subset of a set of vertices of some 
regular $m$-gon contained in $Q$, where $m \leq |R|$.  
The bound $\frac{3}{2}n$ on $|R|$ 
in the statement of the theorem is best possible.
\end{theorem}

We remark that by Observation \ref{observation:directions}, 
if $Q$ is a circle and each of $A$ and $B$ is 
a subset of the set of vertices of a regular $m$-gon contained in $Q$, then
the number of distinct directions of lines passing through a point
in $A$ and a point in $B$ is at most $m$. 
This is equivalent to saying that one can find a set $R$ of
$m$ points on the line at infinity such that every line through a point in $A$
and a point in $B$ passes through a point in $R$. 

\begin{figure}[ht]
   \centering
    \includegraphics[width=7cm]{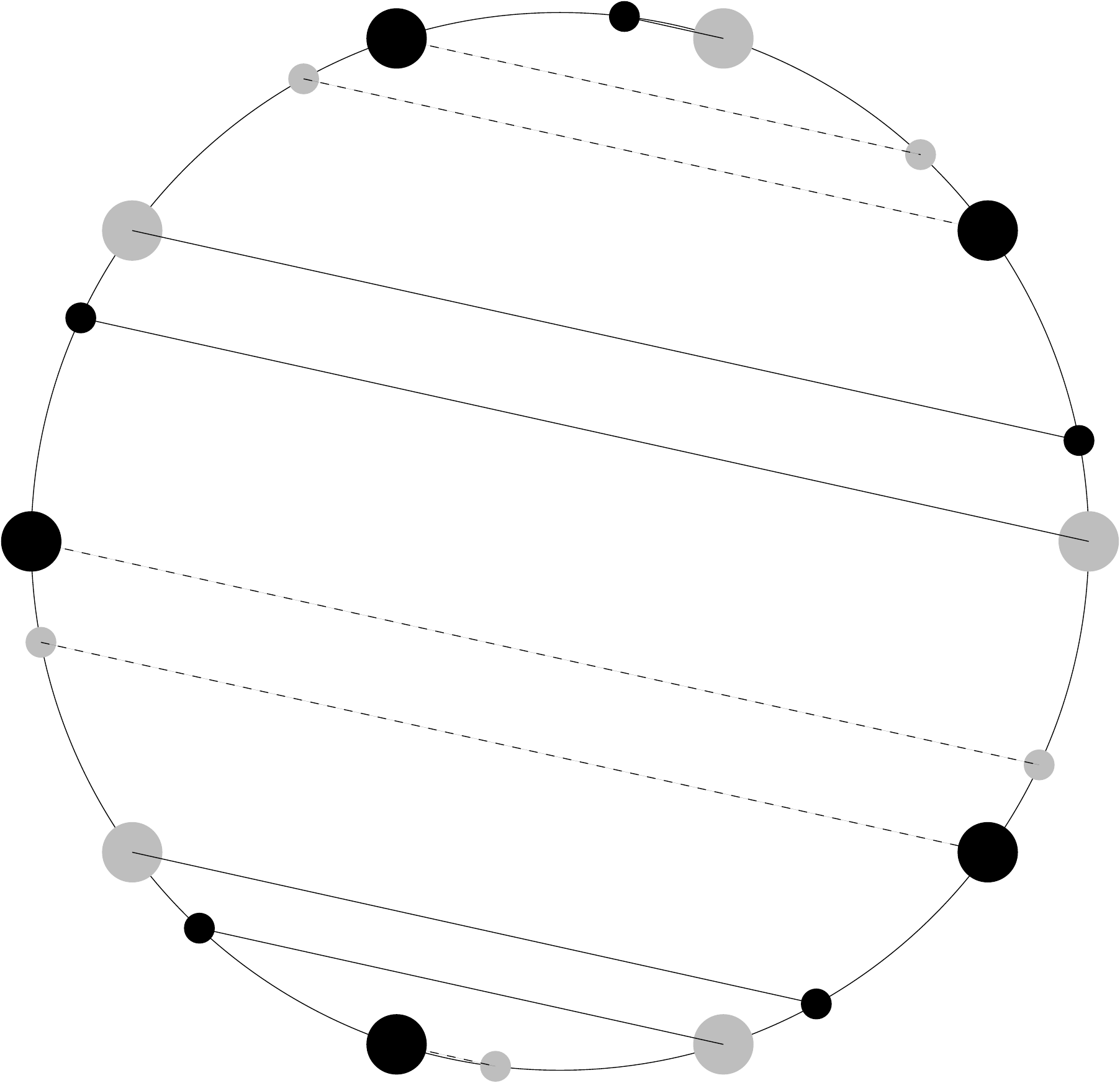}
        \caption{The construction with $|B|=|G|=10$. The points of $B=Z \cup (-Z')$ are 
colored black. The points of $Z$ are drawn by smaller black discs and the points
of $-Z'$ by bigger black discs. 
The points of $G=(-Z) \cup Z'$ are colored gray. The points of $-Z$ are drawn by smaller gray discs and the points
of $Z'$ by bigger gray discs.}
        \label{figure:mainb_reducible}
\end{figure}

Another important remark is about the tightness of the bound 
$|R| < \frac{3}{2}n$ in the statement Theorem \ref{theorem:mainb_reducible}.
We now present the construction showing this.
Let $Q$ be a circle centered at the origin and let $k$ be an odd integer.
Let $Z$ be the set of vertices of a regular $k$-gon contained in~$Q$.
Let $Z'$ be a generic rotation of $Z$ about the center of~$Q$.
Then set $B=Z \cup (-Z')$ and $G=-Z \cup Z'$. 
We have $|B|=|G|=2k$ and set $n=2k$.
We claim that the set of lines passing through a point in 
$B$ and a point in $G$ may appear in one of at most $3k$ distinct directions.
Once this is verified, let~$R$ be the set of $3k$ points on the line
at infinity that correspond to the $3k$ distinct directions of lines
passing through a point in $B$ and a point in $G$. Then every line
passing through a point of $B$ and a point of $G$ will pass also through 
a point of $R$. We have $|R|=\frac{3}{2}n$ and neither $B$ or $G$
is contained in the set of vertices of a regular $m$-gon.
To see that indeed the number of distinct directions of lines passing through 
a point in $B$ and a point in $G$ is equal to $3k=\frac{3}{2}n$ we use
Observation \ref{observation:directions}.
Notice that each of $Z, Z', -Z$, and $-Z'$ is the set of vertices
of a regular $k$-gon inscribed in $Q$. Recall that 
$B=Z \cup (-Z')$ and $G=(-Z) \cup Z'$.
By Observation \ref{observation:directions}, the lines passing through 
a point of $Z$ and a point of $-Z$ have precisely $k$ distinct directions.
The same is true for the lines passing through a point of $-Z'$ and a point of 
$Z'$. The crucial observation is that the set of directions of 
lines passing through a point of $-Z'$ and a point of $-Z$ (and there are 
precisely $k$ such distinct directions)
is precisely the same set of directions of the lines passing through 
a point of $Z'$ and a point of $Z$. This implies that the
number of distinct directions of lines passing through a point of 
$B=Z \cup (-Z')$ and a point of $G=(-Z) \cup Z'$ is equal to $3k=\frac{3}{2}n$,
as desired (see Figure \ref{figure:mainb_reducible}).

\medskip
 
The case of Theorem \ref{theorem:mainb}
is simpler than 
the one of Theorem \ref{theorem:main}. For this reason we start with the
proof of Theorem \ref{theorem:mainb} and Theorem \ref{theorem:mainb_reducible} 
in Section \ref{section:mainb}
and continue to the proof of Theorem \ref{theorem:main} and Theorem 
\ref{theorem:main_reducible} in Section \ref{section:main}.

\section{Proof of Theorem \ref{theorem:mainb} and Theorem 
\ref{theorem:mainb_reducible}.}\label{section:mainb}

We will rely on the following 
result about subsets of abelian groups. 
This is one of several variation of, by now a classical,
result of Frieman \cite{F73} who proved a similar result for general groups.
The Lemma as we bring it below is from \cite{BP18} (see Proposition 2.1 there). 

\begin{lemma}[\cite{BP18}]\label{lemma:A+B}
Suppose $A_{1}$ and $A_{2}$ are two finite subsets of an abelian group $F$.
Let $H$ denote the subgroup of $F$ that is the stabilizer of $A_{1}+A_{2}$.
That is, $H=\{x \in F \mid x+A_{1}+A_{2}=A_{1}+A_{2}\}$.
If $|A_{1}| \geq |A_{2}|$, $|A_{2}| \geq \frac{3}{4}|A_{1}|$, and 
if $|A_{1}+A_{2}| < \frac{3}{2}|A_{1}|$, then $A_{1}+A_{2}$ is equal to a 
coset of $H$. Consequently, each of $A_{1}$ and $A_{2}$ is contained in a coset 
of $H$. 
\end{lemma}

Lemma \ref{lemma:A+B} is stated in \cite{BP18} for \emph{finite} abelian 
groups although only the finiteness of $A_{1}$ and $A_{2}$ is required for the
proof. Moreover, the statement of Lemma \ref{lemma:A+B} in \cite{BP18}
only says that such a subgroup $H$ exists but in the proof $H$ is 
taken to be the stabilizer of $A_{1}+A_{2}$. We will not make use of the
fact that $H$ is the stabilizer of $A_{1}+A_{2}$. The easy proof of Lemma
\ref{lemma:A+B} in \cite{BP18}
relies on Kneser's theorem (see \cite{M65}).

We start with the proof of Theorem \ref{theorem:mainb}, where we assume that the
cubic curve $c$ is irreducible.
In this case we have the group action that is naturally defined on $c$,
where three points on $c$ are collinear if and only if their sum is
equal to $0$.
Because every line through a point in $B$ and a point in $G$ passes
through a point in $R$ we conclude that $R \supset -(B+G)$.

Let $H$ be the subgroup of $c$ that is the stabilizer of $B+G$.
By Lemma \ref{lemma:A+B}, where $B$ and $G$ are in the role of $A_{1}$ and 
$A_{2}$, both $B$ and $G$ are subsets of a coset of $H$.
Moreover, $B+G$ is equal to a coset of $H$. In particular, $|H| \leq |R|$.

If $|R|=n$ it follows that $B+G=-R$. Consequently, $R$ is a coset of $H$
and because $|B|=|G|=|H|=n$ both $B$ and $G$ must be equal to cosets of $H$.
This completes the proof of Theorem \ref{theorem:mainb}

We remark that the case $|R|=n$ in Theorem \ref{theorem:mainb}
is extremely easy even without using Lemma \ref{lemma:A+B}. We give here
the easy argument because it is very short.
Choose $b \in B$ and $g \in G$ and let $B'=-B+b$ and $G'=-G+g$.
Then $0 \in B', G'$ and we have $R+(b+g)=B'+G'$.
Because $|B'|=|G'|=|B'+G'|=n$ and $B',G' \subset B'+G'$
(this is because $0 \in B', G'$), it follows that $B'=G'=B'+G'$ and we 
denote this set by $H$.
$H$ must be a subgroup because $H+H=H$ and $0 \in H$. 
Now the result follows because $B=-B'+b=-H+b=H+b$, 
$G=-G'+g=-H+g=H+g$, and $R=B'-(b+g)=H+b-(b+g)=H-g$.
This completes the proof of Theorem \ref{theorem:mainb} in the case $|R|=n$.
\bbox

\bigskip

We now move on to the proof of Theorem \ref{theorem:mainb_reducible}.
We will need the following lemma that can be found in \cite{GT13}.

\begin{lemma}[Proposition $7.3$ in \cite{GT13}]\label{lemma:GT}
Let $c$ be a cubic curve that
is a union of a quadric $Q$ and a line $\ell$. Then there
is an abelian group $F$ and two mapping $\phi_{Q}:F \rightarrow Q$
and $\phi_{\ell}:F \rightarrow \ell$ such that for $x,y,z \in F$
$x+y+z=0$ if and only if $\phi_{Q}(x), \phi_{Q}(y)$, and $\phi_{\ell}(z)$ are
collinear. Moreover, if $Q$ is a hyperbola, then $F$ is 
isomorphic to $\mathbb{Z}/2\mathbb{Z}\times \mathbb{R}$, if 
$Q$ is a parabola, then $F$ is isomorphic to $(\mathbb{R}, +)$, and 
if $Q$ is an ellipse, then $F$ is isomorphic to $(\mathbb{R}/\mathbb{Z}, +)$ .
\end{lemma}

In the case of Theorem \ref{theorem:mainb_reducible} 
the cubic curve $c$ is a union of a quadric $Q$ and a line $\ell$.
We claim that no point of $B \cup G$ may lie on $\ell$.
To see this, assume to the contrary and without loss of generality that 
$b \in B$ lies on $\ell$. Then by considering the $n$ 
lines through $b$ and the $n$
points in $G$ we conclude that there must be at least $n-1$ points of $R$ 
not on $\ell$ and that means there are at least $n-1$ points of $R$ on $Q$.
On the other hand because both $B$ and $G$ are in general position,
there cannot be more than two points of $B$ and two points of $G$ on $\ell$.
Therefore, there are at least $n-2$ points of $B$ on $Q$ and at least $n-2$
points of $G$ on $Q$. Every line through a point of $B$ on $Q$ and a point
of $G$ on $Q$ must contain a point of $R$ on $\ell$. This implies
at least $n-2$ points of $R$ on $\ell$. This is a contradiction as we assume
$|R|< \frac{3}{2}n < 2n-3$. 

Having shown that $B \cup G \subset Q$ it follows that essentially
$R \subset \ell$. This is because no point on $Q$ may be collinear with 
a point of $B$ and a point of $G$ (both also on $Q$).

We now use Lemma \ref{lemma:GT}. We conclude that there
is an abelian group $F$ and two mapping $\phi_{Q}:F \rightarrow Q$
and $\phi_{\ell}:F \rightarrow \ell$ such that for $x,y,z \in F$
$x+y+z=0$ if and only if $\phi_{Q}(x), \phi_{Q}(y)$, and $\phi_{\ell}(z)$ are
collinear. We know moreover, if $Q$ is a hyperbola, then $F$ is 
isomorphic to $\mathbb{Z}/2\mathbb{Z}\times \mathbb{R}$, if 
$Q$ is a parabola, then $F$ is isomorphic to $(\mathbb{R}, +)$, and 
if $Q$ is an ellipse, then $F$ is isomorphic to $(\mathbb{R}/\mathbb{Z}, +)$.
We know that every point of $B$ (on $Q$) and a point of $G$ (on $Q$)
are collinear with a point of $R$.
Therefore, taking $\tilde{B}=\phi_{Q}^{-1}(B)$, $\tilde{G}=\phi_{Q}^{-1}(G)$,
and $\tilde{R}=\phi_{\ell}^{-1}(R)$, we have 
$\tilde{B} + \tilde{G} \subset -\tilde{R}$ and consequently 
$|\tilde{B} + \tilde{G}| \leq \frac{3}{2}n$.

We apply Lemma \ref{lemma:A+B} with the abelian group $F$ taking
$A_1$ and $A_2$ in Lemma \ref{lemma:A+B} to be $A_1=\tilde{B}$ and $A_2=\tilde{G}$,
respectively. 
Lemma \ref{lemma:A+B} implies (because we have 
$|\tilde{B} + \tilde{G}| < \frac{3}{2}n$) that 
$\tilde{B}+\tilde{G}$
is equal to a coset of a subgroup $H$ of $F$ and each of  
$\tilde{B}$ and $\tilde{G}$ is contained in a 
coset of~$H$.
Notice that the cardinality of $H$ is smaller than or equal to $|R|$ and 
consequently, if $|R|=n$, then each of $\tilde{B}$ and $\tilde{G}$ are 
in fact \emph{equal} to a coset of $H$.

We can now further continue and give a more concrete description of 
$B$ and $G$. Project $\ell$ to the line at infinity. 
The points of $R$ correspond to a collection of $|R|$ 
distinct directions and every line through a point of $B$ and a point of $G$
has one of these directions.

We claim that $Q$ cannot be a parabola or 
a hyperbola. To see this, recall that by Lemma~\ref{lemma:GT}, 
if $Q$ is a hyperbola or a parabola, then the abelian group $F$ is either
$\mathbb{Z}/2\mathbb{Z}\times \mathbb{R}$, or $(\mathbb{R}, +)$, 
respectively. In either case
it cannot have a finite subgroup $H$ of size greater than $2$.

Having shown that $Q$ cannot be a parabola or a hyperbola we conclude that
$Q$ must be an ellipse. In this case the abelian group $F$ is isomorphic to
$(\mathbb{R}/\mathbb{Z}, +)$ and the only finite subgroups it has are
isomorphic to $\mathbb{Z}_{k}$. By applying an affine transformation we may 
assume that $Q$ is a circle. 
We may assume without loss of generality that $0 \in H$
and therefore $H$ is a finite subgroup of $F$. Hence $H$
is isomorphic to $(\mathbb{Z}_{k}, +)$ for some 
$k < \frac{3}{2}n$. We know that $|\tilde{B}|=|\tilde{G}|=n$
and $\tilde{B}+\tilde{G}=H$. 

We claim that $\phi_{Q}(H)$ must be the set of vertices of a regular $k$-gon.
To see this observe that $H$ is a finite subgroup of $F$ that is isomorphic to 
$(\mathbb{Z}_{k}, +)$. From Lemma \ref{lemma:GT} 
we know that if $x+y+z=0$ in $F$, then the line through $\phi_{Q}(x)$ and 
$\phi_{Q}(y)$ is parallel to the direction $\phi_{\ell}(z)$ on the line $\ell$
at infinity. Therefore, the lines through two points of 
$\phi_{Q}(H)$ can be only in $|H|=k$ distinct
directions. It now follows from a well known result of Jamison 
(Theorem 2 in \cite{J86}),
and also very easy to show directly because we know that $\phi_{Q}(H)$ is
contained in a circle, that $\phi_{Q}(H)$ is equal to the set of vertices of 
a regular $k$-gon.
Because each of $\tilde{B}$ and $\tilde{G}$ is a subset of a coset of $H$
we conclude that each of $B$ and $G$ is a rotation about the center of $Q$ 
of some (could be different for $B$ and for $G$) subset of size $n$
of the set of vertices of a regular $k$-gon, where $k < \frac{3}{2}n$.
This concludes the proof of Theorem \ref{theorem:mainb_reducible}.
In the remarks following the statement of Theorem~\ref{theorem:mainb_reducible}
it is shown why the bound $|R| < \frac{3}{2}n$ in the statement of the theorem
cannot be improved even by one unit.
\bbox

\section{Proof of Theorem \ref{theorem:main} and Theorem \ref{theorem:main_reducible}}\label{section:main}

The proof of Theorem \ref{theorem:main} is a bit more involved than 
the proof of Theorem \ref{theorem:mainb}.
The reason is that from the fact that a line through two
points in $P$ passes through a point in $R$ we cannot conclude that
$P+P \subset -R$. This is because we have no information about the
sum of a point in $P$ with itself.

For a subset $A$ of a group we denote by
$A \dsum A$ the set $\{a+a' \mid a, a' \in A, ~~a \neq a'\}$.
Therefore, the conditions in Theorem \ref{theorem:main} imply only
$P \dsum P \subset -R$ rather than $P+P \subset -R$. This difference
makes the proof of Theorem \ref{theorem:main} a bit more challenging.

We will need the following result from \cite{L00}
about subsets of abelian groups.

\begin{theorem}[Theorem 3 in \cite{L00}]\label{theorem:seva}
Suppose $F$ is a group, $A \subset F$, and 
$|A \dot{+} A| \leq \frac{1+\sqrt{5}}{2}|A|-(L+2)$, where $L$, called the 
doubling constant of $F$, is the maximum number of solutions to $x+x=a$
over all $a \in F$.
Then $A \dot{+} A = A+A$.
\end{theorem}

We recall that $P$ is a set of $n$ points in general position in the plane.
Let $R$ be another set of $n$ points, disjoint from $P$ such that 
any line through two points of $P$ passes through a point in $R$.
We assume that $P \cup R$ is contained in a cubic curve $c$ in the plane.
We need to show there is a subgroup $H$ of $c$ of size at most $|R|$ 
such that $P \dot{+} P$ is equal to a coset of $H$.
We have $P \dot{+} P \subset -R$ and consequently 
$P \dot{+} P < \frac{3}{2}n$. 
We apply Theorem 
\ref{theorem:seva}, where the abelian group $F$ is the cubic curve $c$
with its prescribed group structure. 
We notice that the doubling constant of the group $F$ is not greater than $6$.
This is because a solution to $x+x=a$ for $a \in F$ means that there is a line 
through the point $a$ that touches $c$ at $x$ with multiplicity $2$. For an irreducible cubic curve there cannot
be more than $6$ tangent lines through any given point $a$. To prove 
this elementary fact one can project $a$ to infinity and assume without loss of generality that
all the lines through $a$ are horizontal. Then at any point $(x,y)$ 
on $c$ in which 
the tangent line is horizontal we must have 
$\frac{\partial c}{\partial x}(x,y)=0$. However, 
$\frac{\partial c}{\partial x}(x,y)=0$ is a quadric and, by Bezout theorem, it 
intersects the cubic $c$ in at most $6$ points (because $c$ is irreducible).

Having verified that the doubling constant of $F$ is at most $6$, Theorem 
\ref{theorem:seva} implies (because we have $|P \dot{+} P| < 
\frac{3}{2}n < \frac{1+\sqrt{5}}{2}n-8$, where we assume $n$ is large enough) 
that $P \dot{+} P=P+P$.
Now we apply Lemma \ref{lemma:A+B}, where we take $A_{1}=A_{2}=P$ in 
Lemma \ref{lemma:A+B}. We have $|P+P|=|P \dot{+} P| < \frac{3}{2}|P|$. 
Therefore, by Lemma \ref{lemma:A+B} $P \dot{+} P=P+P$
is equal to a coset of a subgroup $H$ and $P$ is contained in a coset of $H$.
Notice that the cardinality of $H$ is smaller than or equal to $|R|$ and consequently, if $|R|=n$, then $P$ is equal to a coset of $H$.
\bbox

\bigskip

We now move on to the proof of Theorem \ref{theorem:main_reducible}.
We start by following the proof of Theorem~\ref{theorem:mainb_reducible}.
The cubic curve $c$ is a union of a quadric $Q$ and a line $\ell$.
We claim that no point of $P$ may lie on $\ell$.
To see this, assume to the contrary that 
$x \in P$ lies on $\ell$. Then by considering the $n-1$ 
lines through $x$ and the $n-1$ other points in $P$
we conclude that there must be at least $n-2$ points of $R$ 
not on $\ell$ and that means there are at least $n-2$ points of $R$ on $Q$.
On the other hand, because $P$ is in general position,
there cannot be more than two points of $P$ on $\ell$.
Therefore, there are at least $n-2$ points of $P$ on $Q$. 
Every line through two points of $P$ on $Q$ must contain a point of $R$ on 
$\ell$. This implies
at least $n-3$ points of $R$ on $\ell$. This is a contradiction as we assume
$|R| < \frac{3}{2}n < (n-2)+(n-3)$.

Having shown that $P \subset Q$ it follows that essentially
$R \subset \ell$. This is because a line through two points (of $P$) on $Q$
cannot be collinear with another point (of $R$) on $Q$.

We now use Lemma \ref{lemma:GT}. We conclude that there
is an abelian group $F$ and two mappings $\phi_{Q}:F \rightarrow Q$
and $\phi_{\ell}:F \rightarrow \ell$ such that for $x,y,z \in F$ we have 
$x+y+z=0$ if and only if $\phi_{Q}(x), \phi_{Q}(y)$, and $\phi_{\ell}(z)$ are
collinear. We know moreover, if $Q$ is a hyperbola, then $F$ is 
isomorphic $\mathbb{Z}/2\mathbb{Z}\times \mathbb{R}$, if 
$Q$ is a parabola, then $F$ is isomorphic to $(\mathbb{R}, +)$, and 
if $Q$ is an ellipse, then $F$ is isomorphic to $(\mathbb{R}/\mathbb{Z}, +)$.
We know that every two points of $P$ (on $Q$) are collinear with a point 
of $R$ (on $\ell$). Therefore, taking $\tilde{P}=\phi_{Q}^{-1}(P)$
and $\tilde{R}=\phi_{\ell}^{-1}(R)$, we have 
$\tilde{P} \dot{+} \tilde{P} \subset -\tilde{R}$ and consequently 
$|\tilde{P} \dot{+} \tilde{P} |\leq \frac{1+\sqrt{5}}{2}n-4$.
We can now continue similarly to the proof of Theorem 
\ref{theorem:mainb_reducible}.
We apply Theorem 
\ref{theorem:seva} with the abelian group $F$. 
We notice that the doubling constant of the group~$F$ is not greater than $2$.
Indeed, consider a solution to $x+x=a$ for $a \in F$. This means that 
$x+x+(-a)=0$ but this means that $\phi_{\ell}(-a)$ is a point on $\ell$
and the line through it and $\phi_{Q}(x)$ is tangent to $Q$. As $Q$ is a conic
there cannot be more than two such points $x$.
Having verified that the doubling constant of $F$ is at most $2$, Theorem 
\ref{theorem:seva} implies (because we have $|\tilde{P} \dot{+} \tilde{P}| < 
\frac{3}{2}n \leq 
\frac{1+\sqrt{5}}{2}n-4$, where $n$ is large enough) that 
$\tilde{P} \dot{+} \tilde{P}=\tilde{P}+\tilde{P}$.
Now we apply Lemma~\ref{lemma:A+B}, where we take $A_{1}=A_{2}=\tilde{P}$ in 
Lemma \ref{lemma:A+B}. We have $|\tilde{P}+\tilde{P}|=
|\tilde{P} \dot{+} \tilde{P}| < \frac{3}{2}|\tilde{P}|$. Therefore, by Lemma \ref{lemma:A+B}, 
$\tilde{P} \dot{+} \tilde{P}=\tilde{P}+\tilde{P}$
is equal to a coset of a subgroup $H$ of $F$ and $\tilde{P}$ is contained in a 
coset of $H$.
Notice that the cardinality of $H$ is smaller than or equal to $|R|$ and consequently, if $|R|=n$, then $\tilde{P}$ is equal to a coset of $H$.

We can now further continue and give a more concrete description of 
$P$. Project $\ell$ to the line at infinity. 
The points of $R$ correspond to a collection of $|R|$ 
distinct directions and every line through
two points of $P$ has one of these directions.

We claim that $Q$ cannot be a parabola or 
a hyperbola. To see this, recall that by Lemma~\ref{lemma:GT}, 
if $Q$ is a hyperbola or a parabola, then the abelian group $F$ is either
$\mathbb{Z}/2\mathbb{Z}\times \mathbb{R}$, or $(\mathbb{R}, +)$, 
respectively. In either case
it cannot have a finite subgroup $H$ of size greater than $2$.

Having shown that $Q$ cannot be a parabola or a hyperbola we conclude that
$Q$ must be an ellipse. In this case the abelian group $F$ is isomorphic to
$(\mathbb{R}/\mathbb{Z}, +)$ and the only finite subgroups it has are
isomorphic to $\mathbb{Z}_{k}$. By applying an affine transformation we may 
assume that $Q$ is a circle. 
We may assume without loss of generality that $0 \in \tilde{P}$
and therefore, $H$ is a finite subgroup of $F$. Hence $H$
is isomorphic to $(\mathbb{Z}_{k}, +)$ for some 
$k < \frac{3}{2}n$. We know that $|\tilde{P}|=n$
and $\tilde{P}+\tilde{P}=H$. The geometric consequence is now clear.
$P$ must be a subset of the set of vertices of a regular $k$-gon on $Q$.
To see this observe that $H$ is a finite subgroup of $F$ that is isomorphic to 
$(\mathbb{Z}_{k}, +)$. From Lemma \ref{lemma:GT} 
we know that if $x+y+z=0$ in $F$, then the line through $\phi_{Q}(x)$ and 
$\phi_{Q}(y)$ is parallel to the direction $\phi_{\ell}(z)$ on the line $\ell$
at infinity. Therefore, the lines through two points of 
$\phi_{Q}(H)$ can be only in $|H|$ distinct
directions. As in the proof of Theorem \ref{theorem:mainb_reducible},
it now follows from a well known result of Jamison (\cite{J86}),
and also very easy to show directly, that $\phi_{Q}(H)$ is affinely equivalent 
to the set of vertices of a regular $|H|$-gon. This concludes
the proof of Theorem \ref{theorem:main_reducible}.
In the remarks following the statement of Theorem \ref{theorem:main_reducible}
it is shown why the bound $|R| < \frac{3}{2}n$ in the statement of the theorem
cannot be improved even by one unit.
\bbox

\paragraph{Acknowledgments}
We thank Vsevolod F. Lev for incredibly useful references and up to date 
(and beyond) relevant background about abelian groups.

\bibliographystyle{plain}

\end{document}